\documentclass[a4paper,11pt]{article} 

\usepackage[paper=a4paper,left=1.8cm,right=1.8cm]{geometry}
\usepackage{graphics} 
\usepackage{epsfig} 
\usepackage{times} 
\usepackage{amsmath} 
\usepackage{amssymb} 

\usepackage{epstopdf}
\usepackage{theorem}
\usepackage[]{units}
\usepackage{subfig}
\usepackage{psfrag}
\usepackage{graphicx}
\usepackage{mathabx}
\usepackage{ stmaryrd }
\usepackage{wrapfig}
\usepackage{color}
\usepackage[percent]{overpic}
\usepackage{tabularx}

 \newtheorem{lemma}{Lemma}
 \newtheorem{corollary}{Corollary}
 \newtheorem{proposition}{Proposition}
 
 \newtheorem{example}{Example}

  \newtheorem{proof}{Proof}

\newcommand{\R}{\mathbb{R}}
\newcommand{\polytope}{\Pi}

\title{Accelerating MPC by online detection of state space sets with common optimal feedback laws}
\author{Kai K\"onig and Martin M\"onnigmann\thanks{Corresponding author.}\\
Automatic Control and Systems Theory, Department of Mechanical Engineering,\\
	Ruhr-Universit\"at Bochum, 44801 Bochum, Germany.\\ E-mail: {\tt\small kai.koenig-h4d@rub.de} and {\tt\small martin.moennigmann@rub.de}}

\begin{document}

\maketitle

\begin{abstract}
	Model predictive control (MPC) samples a generally unknown and complicated feedback law point by point. 
	The solution for the current state $x$ contains, however, 
	more information than only the optimal signal $u$ for this particular state. 
	In fact, it provides an optimal affine feedback law $x\rightarrow u(x)$ on a polytope $\Pi \subset\R^n$, 
	i.e., on a full-dimensional state space set. 
	It is an obvious idea to reuse this affine feedback law as long as possible.
	Reusing it on its polytope $\Pi$ is too conservative, however, 
	because any $\Pi$ is a state space set with a common affine law 
	$x\rightarrow (u_0^\prime (x), \dots, u_{N-1}^\prime (x))^\prime\in\R^{Nm}$ 
	for the entire horizon $N$. We show a simple criterion exists for identifying the polytopes that have a common $x\rightarrow u_0(x)$, but may differ with respect to $u_1(x), \dots, u_{N-1}(x)$.
	Because this criterion is
	too computationally expensive for an online use, 
	we introduce a simple heuristics for the fast construction of a subset of the polytopes of interest. 
	Computational examples show (i) a considerable fraction of QPs can be avoided ($\unit[10]{\%}$ to $\unit[40]{\%}$) 
	and (ii) the heuristics results in a reduction very close to the maximum one that could be achieved if the explicit solution was available. 
	We stress the proposed approach is intended for use in online MPC and it does not require the explicit solution. 
\end{abstract}

\section{Introduction}\label{sec:Intro}

Many efforts have been devoted to reducing the computational cost of model predictive control (MPC). 
We do not treat numerical methods such as tailored optimization algorithms here, 
but exploit the piecewise-affine structure of the solution~\cite{Bemporad2002,Tondel2001,Seron2003}. 
We stress that we never calculate explicit control laws, but the present paper belongs to a group of works 
\cite{Kvasnica2017, Li2014, Lehmann2013, Eqtami2011, Jost2013}
that exploit the affine structure, or the corresponding structure of the set of active sets~\cite{Gupta2011,Feller2013,Herceg2015,Oberdieck2017,Monnigmann2019},
to accelerate online MPC.

Specifically, the present paper extends regional MPC~\cite{Jost2015a, Koenig2017a}.
It is the central idea of regional MPC to reuse the optimal solution found by solving the underlying optimal control problem for the current state. 
Regional MPC is based on 
the fact that the pointwise solution of a MPC problem, or more precisely the active set for the current state, defines an affine feedback law and the polytope on which this law provides the optimal solution. Instead of solving an optimization problem in every time step, regional MPC attempts to reuse this feedback law whenever possible. As a consequence, the number of optimization problems to be solved can be reduced without having to know the explicit solution. 
We stress again that we do not precompute feedback laws and polytopes.

Simulations show that closed-loop trajectories often step from polytope to polytope and therefore optimal feedback laws are not reused~\cite{Koenig2017a,Koenig2017c}.  
It is therefore of special interest to extend the region of validity for a feedback law in order to increase the reuse frequency. 
It is well-known that neighboring polytopes often have the same feedback law~\cite{Geyer2008, Kvasnica2012, Kvasnica2013,Oravec2013}.  
This insight has already been exploited for an a posteriori complexity reduction of explicit control laws. In Geyer et al. \cite{Geyer2008}, regions with identical feedback laws are merged based on a hyperplane representation of the explicit solution. In Kvasnica et al. \cite{Kvasnica2012,Kvasnica2013}, regions with a saturated feedback law are eliminated by using simple functions instead. In Oravec et al. \cite{Oravec2013}, regions with the same feedback law are replaced by bounding polygons and their inner and outer approximations. All of these methods require the explicit solution to be known. 

We present an approach for identifying polytopes with the same feedback law that does not require the explicit solution to the optimal control problem. 
We show that some feedback laws are uniquely defined by only a subset of the constraints. Whenever multiple polytopes have this subset in common, they define the same feedback law. 
The computational effort for constructing all polytopes that share the same feedback law with the current one is, however, too high for online use. 
We therefore propose a heuristics for constructing a subset of these polytopes that can be used in online MPC. 

We state the system and problem class along with some preliminaries in Section~\ref{sec:problemstatement}. 
The criterion for the detection of polytopes that have a common optimal MPC feedback law is presented in Section~\ref{sec:ActiveSubsets}, followed by the heuristics for its online use in Section~\ref{sec:Approaches}. We apply the proposed approach to three examples in Section~\ref{sec:example} and give conclusions and an outlook in Section~\ref{sec:conclusion}. 

\section{Problem statement and preliminaries}\label{sec:problemstatement}
\subsection{Optimal control problem (OCP) and piecewise affine structure of the solution}
We consider the optimal control problem   
\begin{subequations}\label{eq:OCP}
	\begin{align}
		&\min \limits_{\substack{u(k),k=0,\ldots,N-1\\x(k),k=1,\ldots,N}}  \quad  \| {x(N)}\|_P^2+\sum \limits_{k=0}^{N-1} \left(
		  \|{x(k)} \|_{Q}^2+\|{u(k)} \|_{R}^2 
		\right)\label{eq:costFunction}\\
		&\text{s.t.} 
		  \ \ \, {x}(k+1)=A{x}(k)+B{u}(k), \ k=0, \ldots, N-1, \label{eq:Sys}\\
		&\quad \ \ \ {x}(k) \in \mathcal{X}, \quad
		{u}(k) \in \mathcal{U}, \quad \quad \ \ k=0, \ldots, N-1,\\
		&\ \ \ \ \ \ \,  {x}(N) \in \mathcal{T}
	\end{align}
\end{subequations}
that is periodically solved for the given initial condition $x(0)$ and horizon $N$ to stabilize the origin of the discrete-time constrained linear system \eqref{eq:Sys} with state variables $x(k) \in \R^n$ and input variables $u(k) \in \R^m$, where the matrices in~\eqref{eq:OCP} have the obvious dimensions. We assume stabilizability of the pair $(A, B)$, detectability of the pair $(Q^{\frac{1}{2}},A)$, $Q \succeq 0$ and $R \succ 0$. Moreover, we assume $\mathcal{X}$, $\mathcal{U}$ and $\mathcal{T} \subseteq \mathcal{X}$ are compact polytopes that contain the origin as an interior point, where a polytope is the intersection of a finite number of halfspaces.
Let $q$ refer to the number of halfspaces in \eqref{eq:OCP}, i.e., inequalities in \eqref{eq:QP} below, and let $q_{\mathcal{U}}$, $q_{\mathcal{X}}$ and $q_\mathcal{T}$ refer to the number of halfspaces required to define $\mathcal{U}$, $\mathcal{X}$ and $\mathcal{T}$, respectively.        
We make the standard choice for $P$ and $\mathcal{T}$ to guarantee asymptotic stability, i.e., $P \succ 0$ is chosen to be the solution of the discrete time Riccati equation and $\mathcal{T}$ is the closure of the largest open set on which the linear quadratic regulator stabilizes the system without activating the constraints. The terminal set $\mathcal{T}$ is calculated with the procedure proposed in Gilbert and Tan \cite{Gilbert1991}. By substituting \eqref{eq:Sys} into the cost function and constraints, problem \eqref{eq:OCP} can be transformed into the equivalent quadratic program (QP)
\begin{align}\label{eq:QP}
		\begin{split}
		&\min \limits_{\bar{U}} \ \frac{1}{2}\bar{U}'H\bar{U}+x(0)'F\bar{U}+\frac{1}{2}x(0)'Yx(0) \\ 
		&\text{s.t.} \quad G\bar{U} \leq w+Ex(0),
	\end{split}
\end{align}
where $\bar{U}=({u}(0)', \ldots, {u}(N-1)')'$ and where the state sequence in \eqref{eq:OCP} can be determined with \eqref{eq:Sys}. A closed-loop system results from solving problem \eqref{eq:QP} in every time step for the current state $x(0)$ and applying the first input of the predicted input sequence, i.e. $u(0)$, to system \eqref{eq:Sys}. We often use $x$ instead of $x(0)$ below for simplicity. 

Let $\mathcal{X}_f$ refer to the set of initial states for which problem \eqref{eq:QP} has a solution. Under the assumptions stated for problem \eqref{eq:OCP}, $H$ is positive definite and there exist a unique optimal input sequence $\bar{U}^\star(x)$ for every $x \in \mathcal{X}_f$. It is known that the optimal solution $\bar{U}^\star: \mathcal{X}_f\rightarrow\R^{mN}$ is a continuous piecewise affine function on a partition of $\mathcal{X}_f$ into a finite number of polytopes $\polytope_1$,  $\polytope_2$, $\ldots$ \cite{Bemporad2002}. 
We call a single affine piece of the piecewise affine function, i.e. 
\begin{align}\label{eq:affineControlLaw}
x \mapsto \bar{K}_j^\star x+ \bar{b}_j^\star \quad \forall x \in \polytope_j
\end{align}
with $\bar{K}^{\star}_j \in \R^{mN \times n}$ and $\bar{b}^{\star}_j \in \R^{mN}$, a \textit{control law}, where we often omit the index $j$ for simplicity.
We call the first $m$ elements of a control law \eqref{eq:affineControlLaw} MPC feedback law, or \textit{feedback law} for short, and denote it
\begin{align}\label{eq:affineFeedbackLaw}
x \mapsto  K^\star x + b^\star \quad \forall x \in \polytope
\end{align}
with $K^\star=\bar{K}^\star_{\{1, \ldots, m\}}$ and $b^\star=\bar{b}^\star_{\{1, \ldots, m\}}$, where a matrix and vector with a set index refer to the obvious submatrix and subvector. 

\subsection{Regional predictive control}\label{sec:regionalMPC}
Regional predictive control as proposed in Jost et al. \cite{Jost2015a} exploits the piecewise affine structure of the optimal solution without calculating it explicitly. 
We summarize the aspects of the approach needed in the present paper. 
Let $\mathcal{A}(x)$ and $\mathcal{I}(x)$, or $\mathcal{A}$ and $\mathcal{I}$ for short, refer to the sets of active and inactive constraints
\begin{align}
\begin{split}
\mathcal{A}(x)&=\{i \in \mathcal{Q}~|~G_i \bar{U}^\star(x)=w_i+E_ix \},\\
\mathcal{I}(x)&=\{i \in \mathcal{Q}~|~G_i\bar{U}^\star(x)<w_i+E_ix \}= \mathcal{Q}\backslash\mathcal{A}(x),\\
\end{split}
\label{eq:Sets}
\end{align}
where $\mathcal{Q}=\{1, \ldots, q \}$ denotes the set of all constraint indices.
We say $\mathcal{A}$ \textit{exists for~\eqref{eq:OCP}} and~\eqref{eq:QP} if there exists an $x\in\mathcal{X}_f$ with active set $\mathcal{A}$. 
We call an $\mathcal{A}\subseteq\mathcal{Q}$ a \textit{candidate} active set, if it is unknown whether it exists or not. 
For any $\mathcal{A}$ and under the assumption that $G_{\mathcal{A}}$ has full row rank, we introduce 
\begin{align}\label{eq:Kb}
\begin{split}
\bar{K}^\star&=H^{-1}(G_{\mathcal{A}})'WS_{\mathcal{A}}-H^{-1}F',\\
    \bar{b}^\star&=H^{-1}(G_{\mathcal{A}})'Ww_{\mathcal{A}},
    \end{split}
\end{align}
and the polytope 
\begin{align}\label{eq:Polytope}
  \polytope(\mathcal{A})=\{ x \in \R^n ~|~T x \leq d \}
\end{align}
with 
    \begin{align}\label{eq:Td}
    \begin{split}
       T&=\begin{pmatrix} WS_{\mathcal{A}} \\
       G_{\mathcal{I}}H^{-1}(G_{\mathcal{A}})'WS_{\mathcal{A}}-S_{\mathcal{I}}\end{pmatrix}, \\
      d&=-\begin{pmatrix} Ww_{\mathcal{A}} \\
      G_{\mathcal{I}}H^{-1}(G_{\mathcal{A}})'Ww_{\mathcal{A}}-w_{\mathcal{I}} \end{pmatrix}
  \end{split}
\end{align}
where $W=(G_{\mathcal{A}}H^{-1}(G_{\mathcal{A}})')^{-1}$ and $S=E+GH^{-1}F'$ with $S \in \R^{q \times n}$. 

The following lemma follows from Theorem 2 in Bemporad et al. \cite{Bemporad2002}.

\begin{lemma}\label{lem:Bemporad}
	Let $x(0) \in \mathcal{X}_f$ be arbitrary with active set $\mathcal{A}(x(0))$.
	Assume the matrix $G_{\mathcal{A}(x(0))}$ has full row rank.
	Then the affine law $\bar{K}^\star x + \bar{b}^\star$ that results for $\mathcal{A}(x(0))$ with~\eqref{eq:Kb} yields the optimal input sequence $\bar{U}^\star$ on the entire polytope 
	$\polytope(\mathcal{A}\left(x(0))\right)$, i.e.,
	$
	  \bar{U}^\star(x)= \bar{K}^\star x + \bar{b}^\star \mbox{ for all } x\in\polytope\left(\mathcal{A}(x(0))\right).
	$
\end{lemma}

Regional predictive control makes use of Lemma \ref{lem:Bemporad} as follows: The optimal solution to problem \eqref{eq:QP} \textit{for a point} $x \in \mathcal{X}_f$ contains more information than the optimal input sequence $\bar{U}^\star(x)$ for this point. 
If $\bar{U}^\star(x)$ is known for an $x\in\mathcal{X}_f$, the active set $\mathcal{A}(x)$ is uniquely defined by~\eqref{eq:Sets}. 
Consequently, an optimal affine feedback law and the polytope \eqref{eq:affineFeedbackLaw} it is optimal on can be calculated with \eqref{eq:Kb} and \eqref{eq:Td}. It is an obvious idea to reuse this feedback law as long as possible and to solve the optimal control problem~\eqref{eq:OCP} only if necessary.   

We claim without giving details that the computational effort for calculating the matrices in \eqref{eq:Kb} and \eqref{eq:Td} is smaller than for solving a QP~\eqref{eq:QP} (see Berner and M\"onnigmann \cite{BernerP2016} for  details). 

As a final preparation we note that typically only a few constraints are active at any optimal solution. 
If, for example, lower and upper bounds on every input and state variable apply in~\eqref{eq:OCP} and~\eqref{eq:QP}, this implies $|\mathcal{Q}|= 2(n+m)N+ q_\mathcal{T}> 2(n+m)N$.
In order for $G_\mathcal{A}$ to have full row rank, $|\mathcal{A}|\le mN$ must hold, since $G$ has $mN$ columns (see~\eqref{eq:QP}),
which implies 
\begin{equation}\label{eq:CardinalityHelper}
  |\mathcal{A}|\le \frac{1}{2} |\mathcal{Q}|-nN- \frac{1}{2} q_\mathcal{T}.
\end{equation}
In typical cases, $|\mathcal{A}|\ll |\mathcal{Q}|$ applies. 

\section{Detection of polytopes with common feedback laws}\label{sec:ActiveSubsets}
\begin{figure*}[]
\footnotesize
\subfloat[The feedback law $x \mapsto K^\star x+b^\star$ is uniquely defined by the active subset $\tilde{\mathcal{A}}=\{ 1\}$. As a result, it not only holds on one polytope, but on all polytopes that result from active sets that are supersets of $\tilde{\mathcal{A}}=\{1\}$. We emphasize that the union of these polytopes is not known to be convex nor connected in general.]{\begin{overpic}[scale=0.67,tics=10]%
{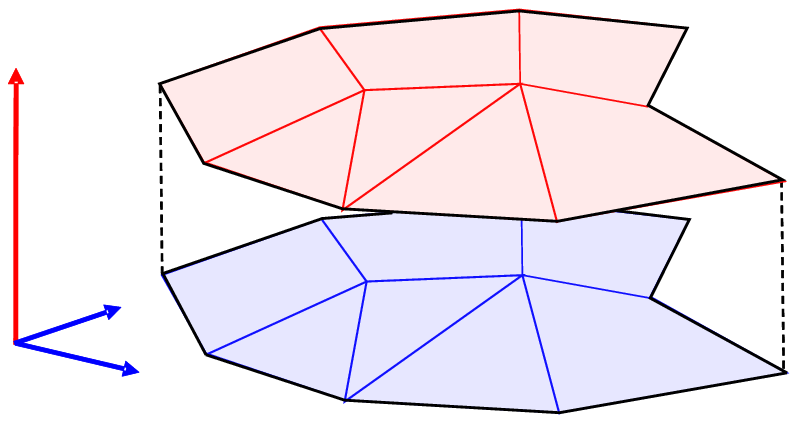}
\put(22,21){$\scriptstyle \{ \textcolor{red}{1},2,3 \}$}
\put(33,12){$\scriptstyle \{ \textcolor{red}{1},2 \}$}
\put(47,17){$\scriptstyle \{ \textcolor{red}{1} \}$}
\put(49,25){$\scriptstyle \{ \textcolor{red}{1},3 \}$}
\put(55,13){$\scriptstyle \{ \textcolor{red}{1},4 \}$}
\put(70,14){$\scriptstyle \{ \textcolor{red}{1},5 \}$}
\put(67,25){$\scriptstyle \{ \textcolor{red}{1},3,4 \}$}
\put(1,50){$ \textcolor{red}{u^{\star}(0)}$}
\put(10,7){$ \textcolor{blue}{x_1}$}
\put(10,23){$ \textcolor{blue}{x_2}$}
\end{overpic}}  \quad
\subfloat[In the existing approach proposed in Jost et al. \cite{Jost2015a}, a QP is solved whenever the current polytope obtained from the pointwise solution has been left (see Sect. \ref{sec:regionalMPC}).
A QP may therefore be solved even if the feedback law does not change. In the sketch, a QP is solved even though the feedback law is constant in all steps.]{\begin{overpic}[scale=0.67,tics=5]%
{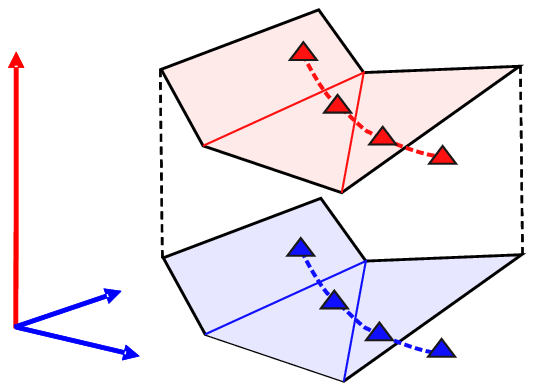}
\put(15,49){$ \textcolor{red}{u^{\star}(0)}$}
\put(26,6){$ \textcolor{blue}{x_1}$}
\put(26,22){$ \textcolor{blue}{x_2}$}
\end{overpic}} \quad
\subfloat[In the approach proposed here (Sect. \ref{sec:Approaches}) the active set $\mathcal{A}=\{ 1,2,3\}$ is calculated by solving a QP. Those subsets of $\mathcal{A}$ which comprise $\tilde{\mathcal{A}}=\{1\}$ yield polytopes, on which the feedback law can be reused. In the sketch, no QP needs to be solved at the points marked by the open circles.]{\begin{overpic}[scale=0.67,tics=10]%
{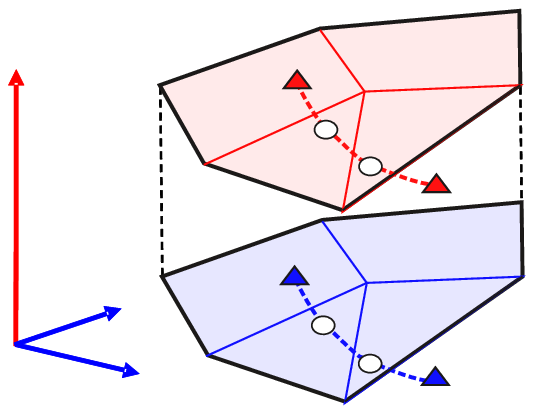}
\put(15,49){$ \textcolor{red}{u^{\star}(0)}$}
\put(39,20){$\textcolor{green}{\scriptstyle \{ 1,2,3 \}}$}
\put(47,14){$\scriptstyle \{ 1,2 \}$}
\put(62,18){$\scriptstyle \{ 1 \}$}
\put(66,26){$\scriptstyle \{ 1,3 \}$}
\put(29,8){$ \textcolor{blue}{x_1}$}
\put(29,22){$\textcolor{blue}{x_2}$}
\end{overpic}}
\captionsetup{width=1\linewidth}
\caption{Comparison of the approach proposed here to Jost et al. \cite{Jost2015a}. Filled triangles denote time steps in which a QP must be solved. White circles denote time steps in which a feedback law from the previous time step can be reused.} \label{fig:idea} 
\end{figure*}

We assume without restriction the constraints in \eqref{eq:OCP} and~\eqref{eq:QP} are ordered such that those on $u(0)$ and $x(1)$ appear first. This can be accomplished, for example, with the order
  \begin{align}\label{eq:ConstraintOrder}
  \begin{split}
  u(0)\in \mathcal{U}, \ x(1) &\in \mathcal{X}, \\
  u(1)\in \mathcal{U}, \ x(2)&\in\mathcal{X}, \\ 
 \vdots \\
 u(N-1)\in\mathcal{U}, \ x(N) &\in \mathcal{T}, \\
 x(0)&\in \mathcal{X}.
 \end{split}
  \end{align}
Essentially, we show that some of the constraints of~\eqref{eq:OCP} and~\eqref{eq:QP} only depend on $u(0)$ but not the remaining $u(k)$ (see Lemma~\ref{lemma:MMOSFirstHelper}) and that this subset of constraints sometimes already determines $u(0)$ uniquely (see Proposition~\ref{prop:uUniquelyDetermined}).
Consequently, certain subsets of active constraints always result in the same MPC feedback law, regardless of the activity of other constraints (see Figure \ref{fig:idea} (a) for a sketch). 
\begin{lemma}\label{lemma:MMOSFirstHelper}
  Consider~\eqref{eq:QP} and assume without restriction the constraints  are ordered as in~\eqref{eq:ConstraintOrder}. 
	Let $\hat{U}=({u}(1)^\prime, \dots, {u}(N-1)^\prime)^\prime$. 
	There exist matrices $G^{11}$, $G^{21}$, $G^{22}$ 
	such that the constraints from~\eqref{eq:QP} have the form
	\begin{equation}\label{eq:ConstraintsBlockForm}
	  \begin{pmatrix} G^{11} & 0 \\ G^{21} & G^{22} \end{pmatrix}
	  \begin{pmatrix} u(0) \\ \hat{U} \end{pmatrix}
	  \le 
	  \begin{pmatrix} w^1 \\ w^2 \end{pmatrix} +
	  \begin{pmatrix} E^1 \\ E^2\end{pmatrix} x(0)
	\end{equation}	  
	where the blocks of $G$ have $q_{\mathcal{U}}+q_{\mathcal{X}}$ and $q- q_{\mathcal{U}}-q_{\mathcal{X}}$ rows and $m$ and $(m-1)N$ columns, respectively. 
\end{lemma}
\begin{proof}
  Since $x(1)= Ax(0)+ Bu(0)$, the first two constraints in~\eqref{eq:ConstraintOrder} are
	\begin{equation}\label{eq:constraintUX}
    	u(0) \in \mathcal{U}, \quad Ax(0)+Bu(0) \in \mathcal{X}.
	\end{equation}
	$\mathcal{U}$ and $\mathcal{X}$ are polytopes by assumption, therefore there exist $G^{11}$, $w^1$ and $E^1$ such that~\eqref{eq:constraintUX}
	is equivalent to 
	\begin{align}\label{eq:firstRowsConstraints}
    	G^{11}u(0) \leq w^1+E^1x(0).
	\end{align} 
	Since $\mathcal{U}$ and $\mathcal{X}$ are defined by $q_{\mathcal{U}}$ and $q_{\mathcal{X}}$ halfspaces, respectively, 
	and since $u(0)$ is of dimension $m$, the block $G^{11}$ has $q_{\mathcal{U}}+q_{\mathcal{X}}$ rows and $m$ columns.
	Relation~\eqref{eq:firstRowsConstraints} is equivalent to the first row in~\eqref{eq:ConstraintsBlockForm}. 
	The second row in~\eqref{eq:ConstraintsBlockForm} collects the remaining $q- q_{\mathcal{U}}- q_{\mathcal{X}}$ constraints from~\eqref{eq:ConstraintOrder}. %\hfill $\square$
\end{proof}

Now let $\mathcal{A}$ be an arbitrary active set that exists for~\eqref{eq:OCP} and~\eqref{eq:QP}
and let $q_\text{stage}= q_\mathcal{U}+ q_\mathcal{X}$. 
Then
\newcommand{\qStage}{q_\text{stage}}
\begin{equation}\label{eq:SufficientActiveSubset} 
  \tilde{\mathcal{A}}:= \mathcal{A}\cap\{1, \dots, \qStage\}
\end{equation}
contains the indices that are active in the first rows of~\eqref{eq:ConstraintsBlockForm}, 
or equivalently, in~\eqref{eq:firstRowsConstraints}. 
If these active constraints, i.e.,
\begin{equation}\label{eq:SufficientConstraints}
  G^{11}_{\tilde{\mathcal{A}}} u(0)= w^1_{\tilde{\mathcal{A}}}+ E^1_{\tilde{\mathcal{A}}} x(0)
\end{equation}
already determine the optimal $u(0)$, then all active sets $\mathcal{A}^\prime$ that exist for~\eqref{eq:OCP} with
\begin{equation}\label{eq:SetsinSolution}
  \mathcal{A}^\prime\cap\{1, \dots, \qStage\}
  = \mathcal{A}\cap\{1, \dots, \qStage\}
  = \tilde{\mathcal{A}}
\end{equation}
determine the same optimal $u(0)$, regardless of the constraints  $i> \qStage$ in the second and subsequent rows of~\eqref{eq:ConstraintOrder}. 
The optimal feedback law must therefore be the same for all active sets in the set
\begin{equation}\label{eq:UnionOfActiveSets}
  \mathcal{M}(\mathcal{A}):= \{
    \mathcal{A}^\prime\subseteq\mathcal{Q}\,|\,
    \mathcal{A}^\prime\cap\{1, \dots, \qStage\}= \mathcal{A}\cap\{1, \dots, \qStage\}
  \}.
\end{equation}
Consequently, the same optimal feedback law applies on the union of polytopes 
\begin{align}\label{eq:UnionOfPolytopes}
  \Gamma\left(\mathcal{A}\right)
  :=\bigcup\limits_{
    \mathcal{A}^\prime\in\mathcal{M}(\mathcal{A})
  }
  \polytope(\mathcal{A}^\prime).
\end{align} 
We can therefore look for sets of active sets that have $m$ independent active constraints on $u(0)$ in common. 
It may appear strange at first sight that such sets of active sets exist,  
because they result in the same optimal feedback signal on the one hand, but differ with respect to active constraints. 
However, the active constraints do not only determine the optimal $u(0)$, but the entire optimal sequence $u(0), \dots, u(N-1)$. 
Essentially, two or more different active sets may therefore have a subset $\tilde{\mathcal{A}}$ in common that determines the same optimal $u(0)$, but the active sets may result in different optimal signals $u(1), \dots, u(N-1)$. 
We summarize the simple criterion that results from the explanations given so far in Proposition~\ref{prop:uUniquelyDetermined}. 

\begin{proposition}\label{prop:uUniquelyDetermined}
  Let $\mathcal{A}$ be an arbitrary active set that exists for~\eqref{eq:OCP} and~\eqref{eq:QP}
  and assume without restriction the constraints are ordered as in~\eqref{eq:ConstraintOrder}. 
  Let $\tilde{\mathcal{A}}$ be as in~\eqref{eq:SufficientActiveSubset}
  and $G^{11}$ as in Lemma~\ref{lemma:MMOSFirstHelper}. If $G^{11}_{\tilde{\mathcal{A}}}$ is invertible,
	the feedback law ${K}^\star x + {b}^\star$ with
	\begin{align}\label{eq:simplifiedFeedbackLaw}
	  {K}^\star=(G^{11}_{\tilde{\mathcal{A}}})^{-1}E_{\tilde{\mathcal{A}}}, 
	  \quad
	  {b}^\star=(G^{11}_{\tilde{\mathcal{A}}})^{-1}w_{\tilde{\mathcal{A}}}
	\end{align}
	yields the same optimal input $u^\star(0)$ as the QP~\eqref{eq:QP} for all $x \in \Gamma(\mathcal{A})$,
	where $\Gamma(\mathcal{A})$ is the union of polytopes defined in~\eqref{eq:UnionOfPolytopes}. 
\end{proposition}
\begin{proof} 
  First note that the constraints can be written in the form~\eqref{eq:ConstraintsBlockForm} and $G^{11}$ is well-defined, 
  since the conditions of Lemma~\ref{lemma:MMOSFirstHelper} are fulfilled. Moreover, 
  the constraints with indices $i\in\mathcal{A}$ hold with equality and $\tilde{\mathcal{A}}\subseteq\mathcal{A}$, therefore~\eqref{eq:SufficientConstraints} holds. 
  Since $(G^{11}_{\tilde{\mathcal{A}}})^{-1}$ exists by assumption, \eqref{eq:SufficientConstraints} is equivalent to
	\begin{align}\label{eq:Um}
	  u(0) = (G^{11}_{\tilde{\mathcal{A}}})^{-1}E^1_{\tilde{\mathcal{A}}} x+(G^{11}_{\tilde{\mathcal{A}}})^{-1} w^1_{\tilde{\mathcal{A}}}, 
	\end{align} 
	which proves~\eqref{eq:simplifiedFeedbackLaw}. 
	Now by definition of $\Gamma(\mathcal{A})$ in~\eqref{eq:UnionOfPolytopes}, $\tilde{\mathcal{A}}\subseteq \mathcal{A}$ for all $x\in\Gamma\left(\mathcal{A}\right)$.
	This implies~\eqref{eq:SufficientConstraints} and the feedback law~\eqref{eq:Um} holds for all $x\in\Gamma\left(\mathcal{A}\right)$. %\hfill $\square$ 
\end{proof}	

Proposition~\ref{prop:uUniquelyDetermined} requires $G^{11}_\mathcal{A}$ to be invertible.  
Since we assume full row rank of $G_\mathcal{A}$, the invertibility of $G^{11}_\mathcal{A}$ can be established by merely checking its dimensions. 
\begin{lemma}
  Assume the conditions of Proposition~\ref{prop:uUniquelyDetermined} to hold. 
  If $G_{\mathcal{A}}$ has full row rank, 
	then $G^{11}_{\tilde{\mathcal{A}}}$ is invertible if and only if $|\tilde{\mathcal{A}}|= m$. 
\end{lemma}
\begin{proof}
	According to Lemma~\ref{lemma:MMOSFirstHelper}, $G$ has the block form stated in~\eqref{eq:ConstraintsBlockForm}. This implies
	\begin{equation*}
	  G_{\mathcal{A}}
	  = \begin{pmatrix} G^{11} & 0 \\ G^{21} & G^{22} \end{pmatrix}_{\mathcal{A}} 
	  = \begin{pmatrix} G^{11}_{\tilde{\mathcal{A}}} & 0 \\ G^{21}_{\hat{\mathcal{A}}} & G^{22}_{\hat{\mathcal{A}}} \end{pmatrix},
	\end{equation*}
	where $\hat{\mathcal{A}}= \mathcal{A}\backslash\tilde{\mathcal{A}}$. 
	Since $G_{\mathcal{A}}$ has full row rank and $\tilde{\mathcal{A}}$ and $\hat{\mathcal{A}}$ partition $\mathcal{A}$, 
	$G^{11}_{\tilde{\mathcal{A}}}$ has full row rank, i.e., row rank $|\tilde{\mathcal{A}}|$. 
	Now assume $|\tilde{\mathcal{A}}|= m$, then $G^{11}_{\tilde{\mathcal{A}}}$ has $m$ independent rows. 
	Since it has $m$ columns according to Lemma~\ref{lemma:MMOSFirstHelper}, it is square and has full rank, which implies invertibility. 
	Conversely, if $G^{11}_{\tilde{\mathcal{A}}}$ is invertible, it must have exactly $m$ independent rows, which implies $|\tilde{\mathcal{A}}|=m$. 
\end{proof}

We stress that a set $\Gamma(\mathcal{A})$ is not in general convex and may not be connected, since it is a union of polytopes. 
This is not a restriction, however, for the method proposed here, which will become evident in Sections~\ref{sec:Approaches} and~\ref{sec:example}.  
Secondly, we stress that we do not require $\Gamma(\mathcal{A})$ in the online approach proposed in Section~\ref{sec:Approaches}, but we determine them only for the sake of comparisons.  
Finally, we note that $\Gamma(\mathcal{A})$ also arises as the intersection of \textit{regions of activity} (see Jost and M\"onnigmann \cite{Jost2013} for a detailed discussion). 

\section{Heuristics for constructing sets of polytopes with common feedback laws}\label{sec:Approaches}
Regional MPC aims at reducing the number of QPs to be solved by reusing MPC feedback laws as long as possible. 
It is an obvious idea to reuse the optimal feedback law for $u(0)$ not only on the polytope $\polytope(\mathcal{A})$ it was constructed for as proposed in Jost et al. \cite{Jost2015a}, 
but on the union of polytopes $\Gamma(\mathcal{A})$ or a subset thereof. 

It is not necessary to determine the entire set $\Gamma(\mathcal{A})$ or all active sets $\mathcal{M}(\mathcal{A})$ that define it. 
Any subset of $\mathcal{M}(\mathcal{A})$ is useful, since such a subset yields polytopes on which the current feedback law remains optimal. 
In the example sketched in Figure~\ref{fig:idea}, for example, the region $\Gamma(\mathcal{A})$ comprises the seven polytopes shown in (a), but the subset of the four polytopes shown in (c) already results in a reduction of QPs. 
In fact, the approach proposed in Jost et al. \cite{Jost2015a} corresponds to the smallest possible subset, i.e., the singleton $\{\mathcal{A}\}\subset\mathcal{M}(\mathcal{A})$. 

We claim a useful subset of $\mathcal{M}(\mathcal{A})$ and thus $\Gamma(\mathcal{A})$ can be constructed from the subsets of $\mathcal{A}$. 
This idea results in a simple heuristics, which we explain in the remainder of the present section. Section~\ref{sec:example} then illustrates the reduction that can be achieved with the proposed heuristics. 

Assume $x\in\mathcal{X}_f$ results in an active set $\mathcal{A}$ 
such that Proposition~\ref{prop:uUniquelyDetermined} applies and therefore the constraints 
$\mathcal{A}\cap\{1, \dots, \qStage \}$ 
already define the optimal feedback. 
By its definition~\eqref{eq:UnionOfPolytopes}, the union of polytopes $\Gamma\left(\mathcal{A}\right)$ results from all active sets $\mathcal{A}^\prime$ 
that have at least the active constraints on $u(0)$ and $x(1)$ in common with $\mathcal{A}$, i.e.,
\begin{equation}\label{eq:SubsetHelper1}
  \mathcal{A}^\prime\supseteq \mathcal{A}\cap\{1, \dots, q_{\mathcal{U}}+ q_{\mathcal{X}}\}.
\end{equation}
Because $|\mathcal{Q}|\gg |\mathcal{A}|$ typically holds (see the comments around~\eqref{eq:CardinalityHelper}) it is usually computationally expensive 
to construct all $\mathcal{A}_j \subseteq\mathcal{Q}$ that respect~\eqref{eq:SubsetHelper1}, but it is much less expensive to construct
\begin{equation}\label{eq:SubsetHelper2}
  \mathcal{A}_j\subseteq\mathcal{A}
\end{equation}  
that respect~\eqref{eq:SubsetHelper1}. 
Combining~\eqref{eq:SubsetHelper1} and~\eqref{eq:SubsetHelper2} yields the set of candidate active sets 
\begin{equation}\label{eq:SetsForSubregion}
  \mathcal{E} (\mathcal{A})= \{
  \mathcal{A}^\prime \subseteq\mathcal{Q}\,|\,
  \mathcal{A}\cap\{1, \dots, q_{\mathcal{U}}+ q_{\mathcal{X}}\} \subseteq \mathcal{A}^\prime
  \subseteq\mathcal{A} 
  \}.
\end{equation}

The sets in $\mathcal{E}(A)$ are only candidates, i.e., they may or may not exist for~\eqref{eq:OCP} and~\eqref{eq:QP}. 
If a set $\mathcal{A}^\prime \in \mathcal{E}(\mathcal{A})$ does not exist in the solution of \eqref{eq:OCP}, $\polytope(\mathcal{A}^\prime)=\emptyset$ holds. 
Otherwise, the feedback law defined by $\mathcal{A}\cap\{1, \dots, q_{\mathcal{U}}+q_{\mathcal{X}}\}$ yields the optimal input $u^{\star}(0)$ on the polytope $\polytope(\mathcal{A}^\prime)$. 
We summarize these statements in Corollary \ref{cor:NonOptimalSets}, which is an immediate consequence to Lemma~\ref{lem:Bemporad}, and in Corollary~\ref{cor:SubregionOfActivity}, which follows from Proposition~\ref{prop:uUniquelyDetermined}.

\begin{corollary}\label{cor:NonOptimalSets}
Let $\mathcal{A}$ be an arbitrary active set that does not exist in the solution for \eqref{eq:OCP} and~\eqref{eq:QP} and assume $G_{\mathcal{A}}$ to have full row rank. 
Then the polytope $\polytope(\mathcal{A})$ defined as in~\eqref{eq:Polytope}  is empty.    
\end{corollary}

\begin{proof} 
	First note that the inverse of $G_{\mathcal{A}}H^{-1}G_{\mathcal{A}}'$ exists, since $G_{\mathcal{A}}$ has full row rank and $H \succ 0$ by assumption. 
	Consequently, the polytope $\polytope(\mathcal{A})= \{x\in\R^n\,|\, Tx\le d\}$ with $T$ and $d$ as defined in~\eqref{eq:Td} is still well-defined. 
	Now assume the active set $\mathcal{A}$ does not exist for \eqref{eq:OCP} and $\polytope(\mathcal{A})\ne\emptyset$, and show this leads to a contradiction: 
	Let $x\in\polytope(\mathcal{A})$ be arbitrary. 
	According to Lemma~\ref{lem:Bemporad} the affine law $\bar{K}^\star x+ \bar{b}^\star$ yields the optimal input sequence $\bar{U}^\star$. 
	Substituting $\bar{U}^\star$ into the constraints of~\eqref{eq:QP} results in the active set $\mathcal{A}$, which is the desired contradiction. 
\end{proof}

Finally, Corollary~\ref{cor:SubregionOfActivity} can be stated, which is the basis for the proposed heuristics.
\begin{corollary}\label{cor:SubregionOfActivity}
Assume the conditions of Proposition~\ref{prop:uUniquelyDetermined} to hold and let 
$\mathcal{E}(\mathcal{A})$ be as in~\eqref{eq:SetsForSubregion}. 
Then, for any $\mathcal{A}^\prime\in\mathcal{E}(\mathcal{A})$ with 
$
  |\mathcal{A}^\prime\cap\{1, \dots, q_\mathcal{U}+q_\mathcal{X}\}|= m,
$
	the feedback law ${K}^\star x + {b}^\star$ with
	\begin{align}\label{eq:SubregionOfActivityFeedbackLaw}
	  {K}^\star&=(G^{11}_{\tilde{\mathcal{A}}^\prime})^{-1}E_{\tilde{\mathcal{A}}^\prime}, 
	  \quad
	  {b}^\star=(G^{11}_{\tilde{\mathcal{A}}^\prime})^{-1}w_{\tilde{\mathcal{A}}^\prime},
	  \\
	  \tilde{\mathcal{A}}^\prime&= \mathcal{A}^\prime\cap\{1, \dots, q_\mathcal{U}+ q_\mathcal{X}\}
	\end{align}
	yields the optimal input $u^\star(0)$ for all $x\in\polytope(\mathcal{A}^\prime)$.
\end{corollary}
\begin{proof} 
Let $\mathcal{A}^\prime\in\mathcal{E}(\mathcal{A})$ be arbitrary. 
Since $\mathcal{A}^\prime\subseteq\mathcal{A}$ by definition of $\mathcal{E}(\mathcal{A})$ and since $G_{\mathcal{A}}$ has full row rank by assumption,
$G_{\mathcal{A}^\prime}$ has full row rank. 
Consequently, $G^{11}_{\tilde{\mathcal{A}}^\prime}$, where $G^{11}$ is as in~\eqref{eq:ConstraintsBlockForm}, has full row rank. 
Since $|\tilde{\mathcal{A}}^\prime|= m$, Lemma~\ref{lemma:MMOSFirstHelper} applies and $G^{11}_{\tilde{\mathcal{A}}^\prime}$ is invertible. 
Now recall $\mathcal{A}^\prime$ is only a candidate active set. If it does not exist in the solution to~\eqref{eq:OCP} and~\eqref{eq:QP}, 
then $\polytope(\mathcal{A}^\prime)= \emptyset$ and the claim holds trivially. 
If $\mathcal{A}^\prime$ exists in the solution to~\eqref{eq:OCP} and~\eqref{eq:QP}, then $\polytope(\mathcal{A}^\prime)\in\Gamma(\mathcal{A})$ and the claim follows with Proposition~\ref{prop:uUniquelyDetermined}. 
\end{proof} 

We use Corollary~\ref{cor:SubregionOfActivity} as follows: 
First we determine the active set $\mathcal{A}$ for a given state $x\in \mathcal{X}_f$ by solving QP \eqref{eq:QP}. 
If $G_\mathcal{A}$ has full row rank, we determine $\tilde{\mathcal{A}}= \mathcal{A}\cap\{1, \dots, q_\mathcal{U}+ q_\mathcal{X}\}$.
If $|\tilde{\mathcal{A}}|= m$ holds, we construct $\mathcal{E}(\mathcal{A})$ according to~\eqref{eq:SetsForSubregion} and apply Corollary~\ref{cor:SubregionOfActivity} as long as possible, i.e., we reuse, without solving any QP, the feedback law~\eqref{eq:SubregionOfActivityFeedbackLaw} as long as the system state remains in the polytopes defined by the active sets $\mathcal{E}(\mathcal{A})$. 
If any of the conditions fails to hold, we solve a QP. 
If the conditions hold, but $|\mathcal{E}(\mathcal{A})|$ is too large to be handled efficiently at runtime, the feedback law~\eqref{eq:SubregionOfActivityFeedbackLaw} may still be reused on the original polytope $\Pi(\mathcal{A})$ following Jost et al. \cite{Jost2015a}. 
In our numerical results in Sect. \ref{sec:example} the latter case never occurred and the computation time was always reduced compared to Jost et al. \cite{Jost2015a}.    

\section{Examples}\label{sec:example}
We illustrate the proposed approaches with three examples. We consider a second order system, because its solution can be visualized, a sixth order system that has served as a benchmark example before \cite{Jost2015b}, and an inverted pendulum on a cart that obviously is open-loop unstable. 
\begin{example}\label{example:SISO}
	Consider the single-input-single-output system with the transfer function
	\begin{align}
  	G(s)=\frac{2}{s^2+s+2} \nonumber
	\end{align}
	that is discretized with the 
	sampling time $T_s=\unit[0.1]{s}$. This results in a system of the form \eqref{eq:Sys} with 
	\begin{align}
		A=\begin{pmatrix}  
		0.8955 & -0.1897 \\
		0.0948 & 0.9903
		\end{pmatrix}, \quad
		B=\begin{pmatrix}
		0.0948 \\
		0.0048
		\end{pmatrix}. \nonumber
	\end{align}
	The example is similar to the one in Seron et al. \cite{Seron2003}, but the system must here respect $-3 \leq x_i \leq 3$, $i=1,2$ 
	and $-2 \leq u_1 \leq 2$ 
	and weighting matrices read $Q=\text{diag}(0.01,4)$ and $R=0.01$. 
	We choose the horizon $N= 4$, which results in a QP with $q=32$ inequalities and 4 optimization variables.
\end{example}

\begin{example}\label{example:MIMO}
Consider the multiple-input-multiple-output system with the transfer function
\begin{align}
G(s)=\begin{pmatrix}
\frac{0.05}{36s^2+6s+1} & \frac{0.02(2s+1)}{8s+1}\\
\frac{0.02(2s+1)}{8s+1} & \frac{0.05}{12s^2+3s+1}
\end{pmatrix} \nonumber
\end{align}
that is discretized with the sampling time  $T_s=\unit[1]{s}$  resulting in a system of the form \eqref{eq:Sys} with $n=6$ and $m=2$. The system must satisfy $-15 \leq x_i \leq 15, i=1,\ldots,6$ and $-20 \leq u_j \leq 20, j=1,2$ and weighting matrices are $Q=10 I^{6 \times 6}$ and $R=0.01 I^{2 \times 2}$. We choose the horizon $N=40$, which results in a QP with $q=658$ inequalities and 80 optimization variables.
\end{example}

\begin{example}\label{example:InvertedPendulum}
Consider an inverted pendulum on a cart. The state vector reads $x(k)=(s(k),\varphi(k),\dot{s}(k),\dot{\varphi}(k))'$ with cart position $s(k)$ and pendulum angle $\varphi(k)$. A discretization with the sampling time $T_s=\unit[0.1]{s}$ results in a system of the form \eqref{eq:Sys} with $n=4$, $m=1$ and matrices 
\begin{align}\nonumber
A=\begin{pmatrix}
 1 & -0.0042 & 0.0911 & -0.0001\\
 0 & 1.1084 & 0.0186 & 0.1029\\
 0 & -0.0826 & 0.8265 & -0.0037\\
 0 & 2.1958 & 0.366 & 1.0941
 \end{pmatrix}, \quad 
 B=\begin{pmatrix}
 0.0014\\
 -0.003 \\
 0.0274 \\
 -0.0582 \\
 \end{pmatrix}.
\end{align}
The state and input constraints read $-1 \leq s \leq 1$, $- \frac{\pi}{3} \leq \varphi \leq \frac{\pi}{3}$, $-9 \leq \dot{s} \leq 9$, $- 2 \pi \leq \dot{\varphi} \leq 2 \pi$ and $-10 \leq u \leq 10$. The weighting matrices are set to $Q=I^{4 \times 4}$ and $R=0.01$. We choose the horizon $N=10$, which results in a QP with $q=138$ inequalities and $10$ optimization variables.
\end{example}

We determine the terminal state weighting matrix $P$ and the terminal set $\mathcal{T}$ as explained in Section~\ref{sec:problemstatement} for all examples. 

We examine how often a feedback law can be reused in the approach proposed in Sect. \ref{sec:Approaches} and compare results to both, the existing approach from Jost et al. \cite{Jost2015a}, 
and the optimal reuse that could be achieved if the complete set of polytopes $\Gamma(\mathcal{A})$ defined in~\eqref{eq:UnionOfPolytopes} was always known. 
The sets $\Gamma(\mathcal{A})$ are computed offline for all active sets $\mathcal{A}$ of the given problem in the latter approach. Note this is only done for the purpose of comparisons. It is not required for the proposed approach and it is not in general practical to do so. 

\begin{figure*}[]
\centering
\captionsetup[subfloat]{justification=centering}
{% This file is generated by the MATLAB m-file laprint.m. It can be included
% into LaTeX documents using the packages graphicx, color and psfrag.
% It is accompanied by a postscript file. A sample LaTeX file is:
%    \documentclass{article}\usepackage{graphicx,color,psfrag}
%    \begin{document}\input{Seron2003OfflineA}\end{document}
% See http://www.mathworks.de/matlabcentral/fileexchange/loadFile.do?objectId=4638
% for recent versions of laprint.m.
%
% created by:           LaPrint version 3.16 (13.9.2004)
% created on:           16-Jul-2019 13:50:13
% eps bounding box:     15 cm x 11.25 cm
% comment:              
%
\begin{psfrags}%
\psfragscanon%
\scriptsize%
%
% text strings:
\psfrag{s02}[t][t]{\color[rgb]{0.15,0.15,0.15}\setlength{\tabcolsep}{0pt}\begin{tabular}{c}$x_1$\end{tabular}}%
\psfrag{s03}[b][b]{\color[rgb]{0.15,0.15,0.15}\setlength{\tabcolsep}{0pt}\begin{tabular}{c}$x_2$\end{tabular}}%
%
% axes ticklabel color:
\color[rgb]{0.15,0.15,0.15}%
%
% xticklabels:
\psfrag{x01}[t][t]{-2}%
\psfrag{x02}[t][t]{-1.5}%
\psfrag{x03}[t][t]{-1}%
\psfrag{x04}[t][t]{-0.5}%
\psfrag{x05}[t][t]{0}%
\psfrag{x06}[t][t]{0.5}%
\psfrag{x07}[t][t]{1}%
\psfrag{x08}[t][t]{1.5}%
%
% yticklabels:
\psfrag{v01}[r][r]{}%
\psfrag{v02}[r][r]{-0.8}%
\psfrag{v03}[r][r]{}%
\psfrag{v04}[r][r]{-0.4}%
\psfrag{v05}[r][r]{}%
\psfrag{v06}[r][r]{0}%
\psfrag{v07}[r][r]{}%
\psfrag{v08}[r][r]{0.4}%
\psfrag{v09}[r][r]{}%
\psfrag{v10}[r][r]{0.8}%
\psfrag{v11}[r][r]{}%
%
% Figure:
\includegraphics[height=0.28\textwidth, width=0.3\textwidth]{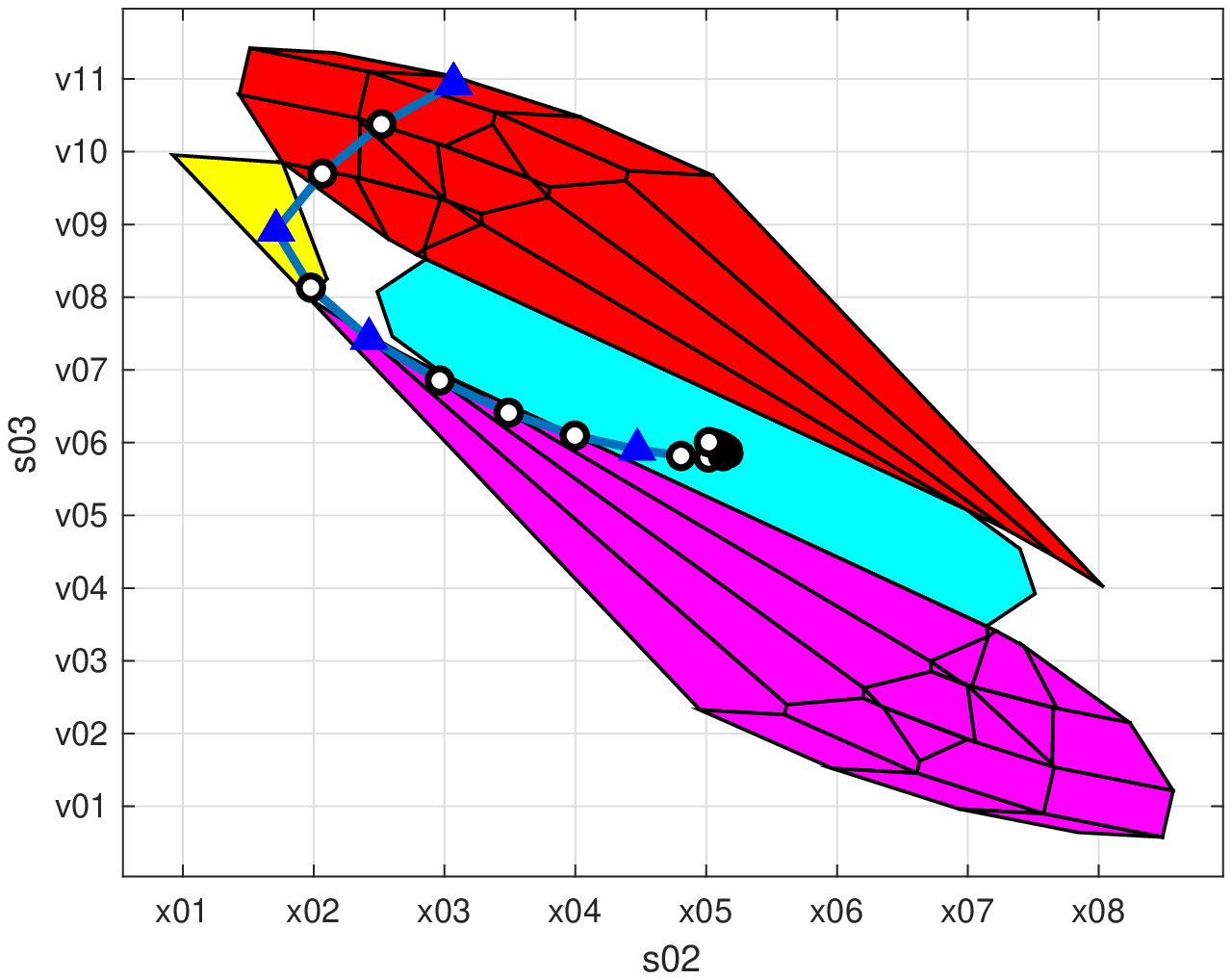}%
\end{psfrags}%
%
% End Seron2003OfflineA.tex
} \quad \quad 
{% This file is generated by the MATLAB m-file laprint.m. It can be included
% into LaTeX documents using the packages graphicx, color and psfrag.
% It is accompanied by a postscript file. A sample LaTeX file is:
%    \documentclass{article}\usepackage{graphicx,color,psfrag}
%    \begin{document}\input{Seron2003OriginalB}\end{document}
% See http://www.mathworks.de/matlabcentral/fileexchange/loadFile.do?objectId=4638
% for recent versions of laprint.m.
%
% created by:           LaPrint version 3.16 (13.9.2004)
% created on:           16-Jul-2019 16:16:39
% eps bounding box:     15 cm x 10.9319 cm
% comment:              
%
\begin{psfrags}%
\psfragscanon%
\scriptsize%
%
% text strings:
\psfrag{s02}[t][t]{\color[rgb]{0.15,0.15,0.15}\setlength{\tabcolsep}{0pt}\begin{tabular}{c}$x_1$\end{tabular}}%
\psfrag{s03}[b][b]{\color[rgb]{0.15,0.15,0.15}\setlength{\tabcolsep}{0pt}\begin{tabular}{c}$x_2$\end{tabular}}%
%
% axes ticklabel color:
\color[rgb]{0.15,0.15,0.15}%
%
% xticklabels:
\psfrag{x01}[t][t]{-2}%
\psfrag{x02}[t][t]{-1.5}%
\psfrag{x03}[t][t]{-1}%
\psfrag{x04}[t][t]{-0.5}%
\psfrag{x05}[t][t]{0}%
\psfrag{x06}[t][t]{0.5}%
\psfrag{x07}[t][t]{1}%
%
% yticklabels:
\psfrag{v01}[r][r]{-0.8}%
\psfrag{v02}[r][r]{}%
\psfrag{v03}[r][r]{-0.4}%
\psfrag{v04}[r][r]{}%
\psfrag{v05}[r][r]{0}%
\psfrag{v06}[r][r]{}%
\psfrag{v07}[r][r]{0.4}%
\psfrag{v08}[r][r]{}%
\psfrag{v09}[r][r]{0.8}%
\psfrag{v10}[r][r]{}%
%
% Figure:
\includegraphics[height=0.28\textwidth, width = 0.3\textwidth]{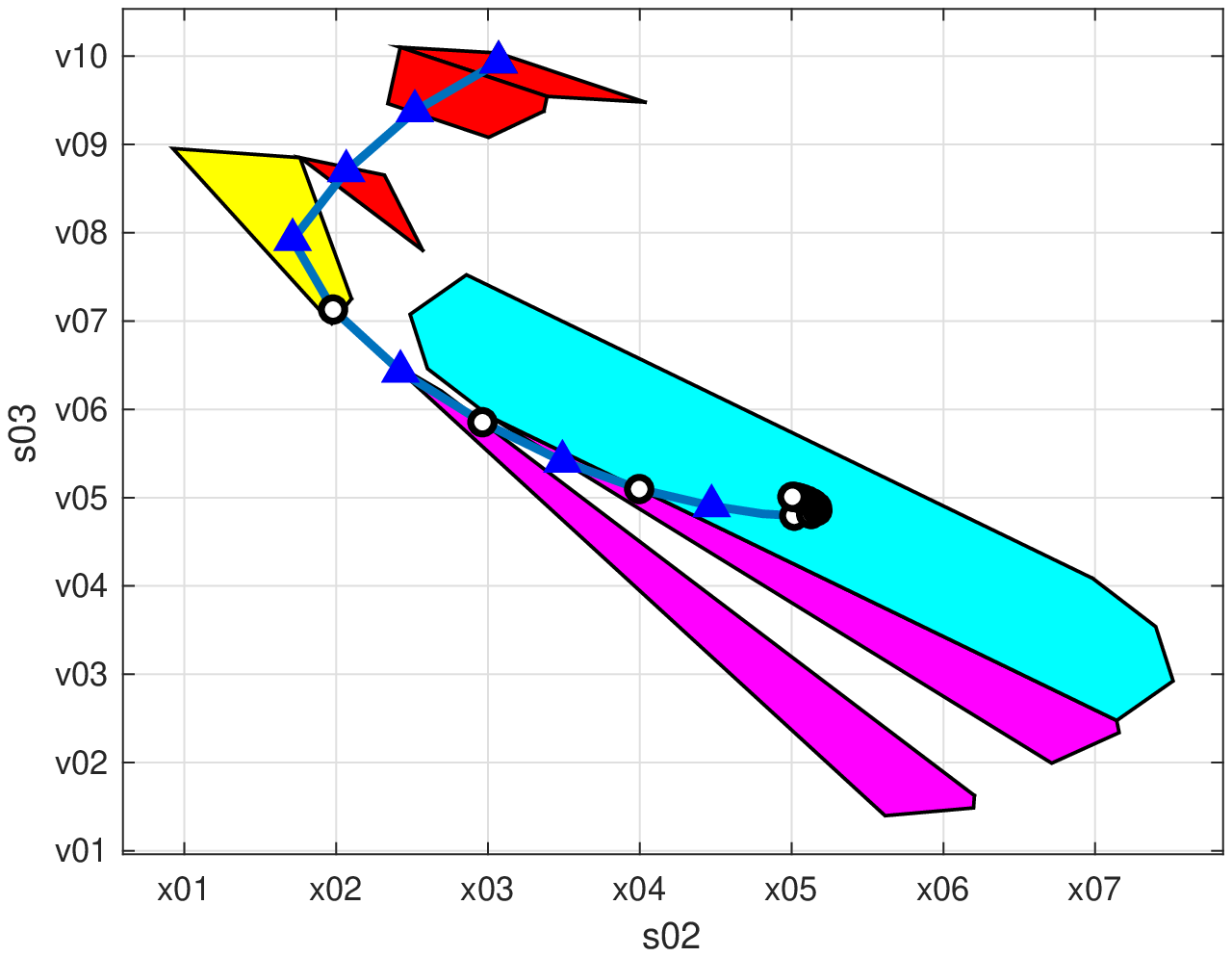}%
\end{psfrags}%
%
% End Seron2003OriginalB.tex
} \quad \quad 
{% This file is generated by the MATLAB m-file laprint.m. It can be included
% into LaTeX documents using the packages graphicx, color and psfrag.
% It is accompanied by a postscript file. A sample LaTeX file is:
%    \documentclass{article}\usepackage{graphicx,color,psfrag}
%    \begin{document}\input{Seron2003OnlineB}\end{document}
% See http://www.mathworks.de/matlabcentral/fileexchange/loadFile.do?objectId=4638
% for recent versions of laprint.m.
%
% created by:           LaPrint version 3.16 (13.9.2004)
% created on:           16-Jul-2019 16:26:22
% eps bounding box:     15 cm x 10.9319 cm
% comment:              
%
\begin{psfrags}%
\psfragscanon%
\scriptsize%
%
% text strings:
\psfrag{s02}[t][t]{\color[rgb]{0.15,0.15,0.15}\setlength{\tabcolsep}{0pt}\begin{tabular}{c}$x_1$\end{tabular}}%
\psfrag{s03}[b][b]{\color[rgb]{0.15,0.15,0.15}\setlength{\tabcolsep}{0pt}\begin{tabular}{c}$x_2$\end{tabular}}%
%
% axes ticklabel color:
\color[rgb]{0.15,0.15,0.15}%
%
% xticklabels:
\psfrag{x01}[t][t]{-2}%
\psfrag{x02}[t][t]{-1.5}%
\psfrag{x03}[t][t]{-1}%
\psfrag{x04}[t][t]{-0.5}%
\psfrag{x05}[t][t]{0}%
\psfrag{x06}[t][t]{0.5}%
\psfrag{x07}[t][t]{1}%
%
% yticklabels:
\psfrag{v01}[r][r]{-0.8}%
\psfrag{v02}[r][r]{}%
\psfrag{v03}[r][r]{-0.4}%
\psfrag{v04}[r][r]{}%
\psfrag{v05}[r][r]{0}%
\psfrag{v06}[r][r]{}%
\psfrag{v07}[r][r]{0.4}%
\psfrag{v08}[r][r]{}%
\psfrag{v09}[r][r]{0.8}%
\psfrag{v10}[r][r]{}%
%
% Figure:
\includegraphics[height=0.28\textwidth, width = 0.3\textwidth]{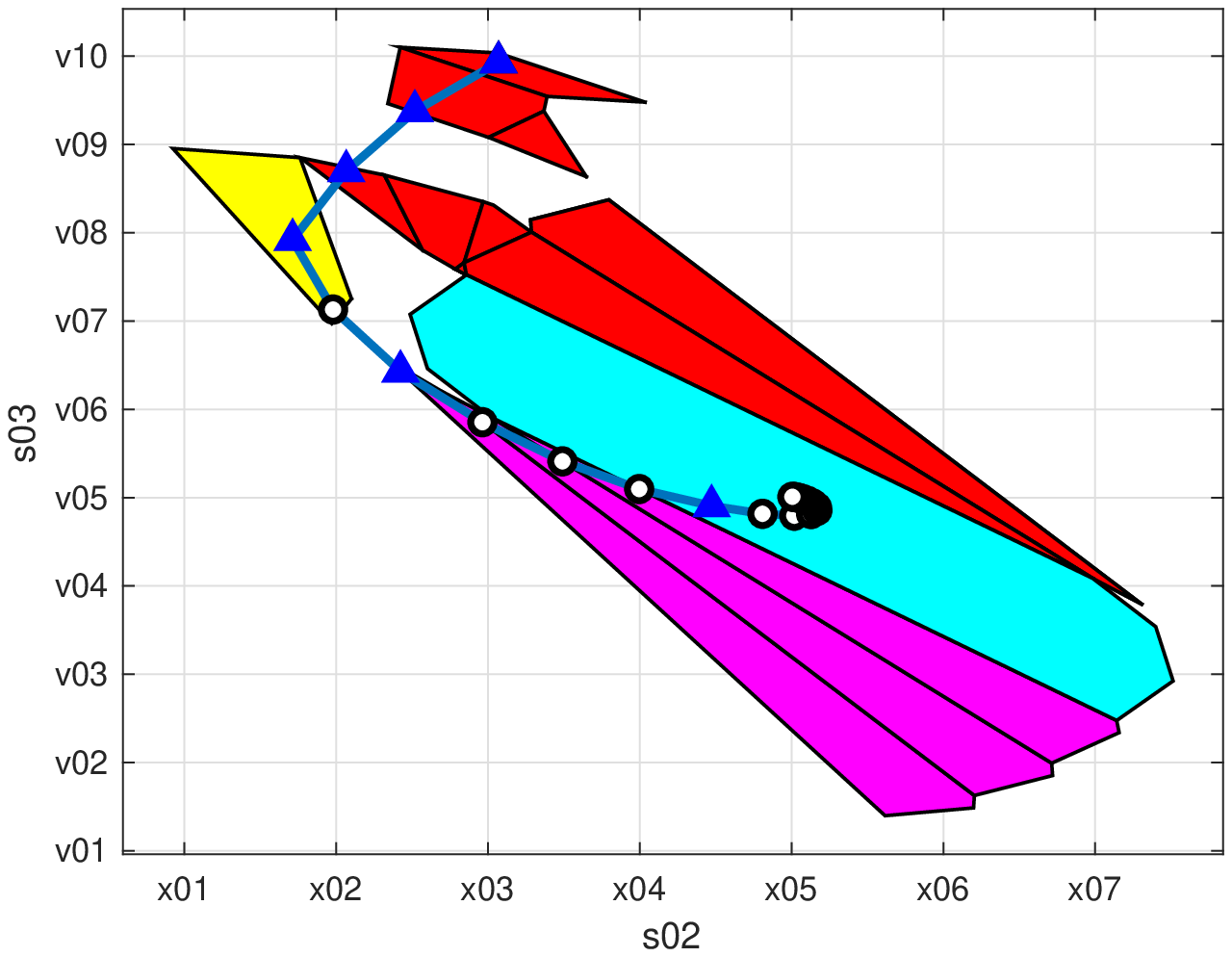}%
\end{psfrags}%
%
% End Seron2003OnlineB.tex
} \\ \vspace{5pt}
\subfloat[known $\Gamma(A)$]{% This file is generated by the MATLAB m-file laprint.m. It can be included
% into LaTeX documents using the packages graphicx, color and psfrag.
% It is accompanied by a postscript file. A sample LaTeX file is:
%    \documentclass{article}\usepackage{graphicx,color,psfrag}
%    \begin{document}\input{Seron2003OfflineB}\end{document}
% See http://www.mathworks.de/matlabcentral/fileexchange/loadFile.do?objectId=4638
% for recent versions of laprint.m.
%
% created by:           LaPrint version 3.16 (13.9.2004)
% created on:           16-Jul-2019 15:59:25
% eps bounding box:     15 cm x 11.25 cm
% comment:              
%
\begin{psfrags}%
\psfragscanon%
\scriptsize%
%
% text strings:
\psfrag{s02}[b][b]{\color[rgb]{0.15,0.15,0.15}\setlength{\tabcolsep}{0pt}\begin{tabular}{c}$x(k)$\end{tabular}}%
\psfrag{s05}[b][b]{\color[rgb]{0.15,0.15,0.15}\setlength{\tabcolsep}{0pt}\begin{tabular}{c}$u(k)$\end{tabular}}%
\psfrag{s07}[t][t]{\color[rgb]{0.15,0.15,0.15}\setlength{\tabcolsep}{0pt}\begin{tabular}{c}$k$\end{tabular}}%
\psfrag{s08}[b][b]{\color[rgb]{0.15,0.15,0.15}\setlength{\tabcolsep}{0pt}\begin{tabular}{c}$e(k)$\end{tabular}}%
%
% axes ticklabel color:
\color[rgb]{0.15,0.15,0.15}%
%
% xticklabels:
\psfrag{x01}[t][t]{0}%
\psfrag{x02}[t][t]{2}%
\psfrag{x03}[t][t]{4}%
\psfrag{x04}[t][t]{6}%
\psfrag{x05}[t][t]{8}%
\psfrag{x06}[t][t]{10}%
\psfrag{x07}[t][t]{12}%
\psfrag{x08}[t][t]{14}%
\psfrag{x09}[t][t]{16}%
\psfrag{x10}[t][t]{18}%
\psfrag{x11}[t][t]{20}%
\psfrag{x12}[t][t]{0}%
\psfrag{x13}[t][t]{2}%
\psfrag{x14}[t][t]{4}%
\psfrag{x15}[t][t]{6}%
\psfrag{x16}[t][t]{8}%
\psfrag{x17}[t][t]{10}%
\psfrag{x18}[t][t]{12}%
\psfrag{x19}[t][t]{14}%
\psfrag{x20}[t][t]{16}%
\psfrag{x21}[t][t]{18}%
\psfrag{x22}[t][t]{20}%
\psfrag{x23}[t][t]{0}%
\psfrag{x24}[t][t]{2}%
\psfrag{x25}[t][t]{4}%
\psfrag{x26}[t][t]{6}%
\psfrag{x27}[t][t]{8}%
\psfrag{x28}[t][t]{10}%
\psfrag{x29}[t][t]{12}%
\psfrag{x30}[t][t]{14}%
\psfrag{x31}[t][t]{16}%
\psfrag{x32}[t][t]{18}%
\psfrag{x33}[t][t]{20}%
%
% yticklabels:
\psfrag{v01}[r][r]{0}%
\psfrag{v02}[r][r]{}%
\psfrag{v03}[r][r]{1}%
\psfrag{v04}[r][r]{-2}%
\psfrag{v05}[r][r]{0}%
\psfrag{v06}[r][r]{2}%
\psfrag{v07}[r][r]{-2}%
\psfrag{v08}[r][r]{0}%
\psfrag{v09}[r][r]{2}%
%
% Figure:
\includegraphics[width = 0.3\textwidth]{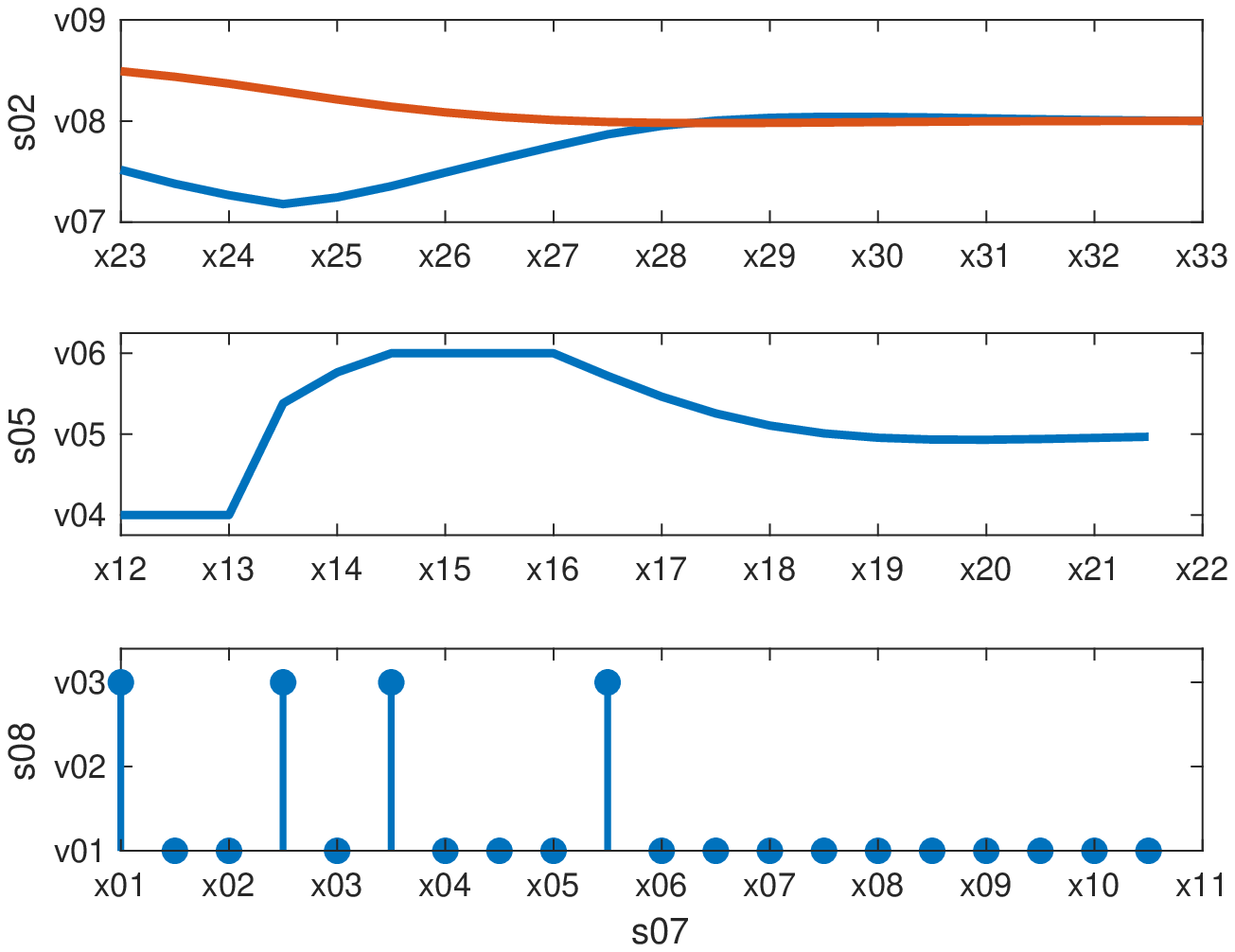}%
\end{psfrags}%
%
% End Seron2003OfflineB.tex
}\quad \quad
\subfloat[existing approach (Jost et al. \cite{Jost2015a})]{% This file is generated by the MATLAB m-file laprint.m. It can be included
% into LaTeX documents using the packages graphicx, color and psfrag.
% It is accompanied by a postscript file. A sample LaTeX file is:
%    \documentclass{article}\usepackage{graphicx,color,psfrag}
%    \begin{document}\input{Seron2003RegA2}\end{document}
% See http://www.mathworks.de/matlabcentral/fileexchange/loadFile.do?objectId=4638
% for recent versions of laprint.m.
%
% created by:           LaPrint version 3.16 (13.9.2004)
% created on:           19-Oct-2018 15:45:17
% eps bounding box:     15 cm x 11.25 cm
% comment:              
%
\begin{psfrags}%
\psfragscanon%
\scriptsize%
%
% text strings:
\psfrag{s03}[b][b]{\color[rgb]{0.15,0.15,0.15}\setlength{\tabcolsep}{0pt}\begin{tabular}{c}$x(k)$\end{tabular}}%
\psfrag{s06}[b][b]{\color[rgb]{0.15,0.15,0.15}\setlength{\tabcolsep}{0pt}\begin{tabular}{c}$u(k)$\end{tabular}}%
\psfrag{s08}[t][t]{\color[rgb]{0.15,0.15,0.15}\setlength{\tabcolsep}{0pt}\begin{tabular}{c}$k$\end{tabular}}%
\psfrag{s09}[b][b]{\color[rgb]{0.15,0.15,0.15}\setlength{\tabcolsep}{0pt}\begin{tabular}{c}$e(k)$\end{tabular}}%
%
% axes ticklabel color:
\color[rgb]{0.15,0.15,0.15}%
%
% xticklabels:
\psfrag{x01}[t][t]{0}%
\psfrag{x02}[t][t]{2}%
\psfrag{x03}[t][t]{4}%
\psfrag{x04}[t][t]{6}%
\psfrag{x05}[t][t]{8}%
\psfrag{x06}[t][t]{10}%
\psfrag{x07}[t][t]{12}%
\psfrag{x08}[t][t]{14}%
\psfrag{x09}[t][t]{16}%
\psfrag{x10}[t][t]{18}%
\psfrag{x11}[t][t]{20}%
\psfrag{x12}[t][t]{}%
\psfrag{x13}[t][t]{}%
\psfrag{x14}[t][t]{}%
\psfrag{x15}[t][t]{}%
\psfrag{x16}[t][t]{}%
\psfrag{x17}[t][t]{}%
\psfrag{x18}[t][t]{}%
\psfrag{x19}[t][t]{}%
\psfrag{x20}[t][t]{}%
\psfrag{x21}[t][t]{}%
\psfrag{x22}[t][t]{}%
\psfrag{x23}[t][t]{}%
\psfrag{x24}[t][t]{}%
\psfrag{x25}[t][t]{}%
\psfrag{x26}[t][t]{}%
\psfrag{x27}[t][t]{}%
\psfrag{x28}[t][t]{}%
\psfrag{x29}[t][t]{}%
\psfrag{x30}[t][t]{}%
\psfrag{x31}[t][t]{}%
\psfrag{x32}[t][t]{}%
\psfrag{x33}[t][t]{}%
%
% yticklabels:
\psfrag{v01}[r][r]{0}%
\psfrag{v02}[r][r]{}%
\psfrag{v03}[r][r]{1}%
\psfrag{v04}[r][r]{-2}%
\psfrag{v05}[r][r]{0}%
\psfrag{v06}[r][r]{2}%
\psfrag{v07}[r][r]{-2}%
\psfrag{v08}[r][r]{0}%
\psfrag{v09}[r][r]{2}%
%
% Figure:
\includegraphics[width=0.3\textwidth]{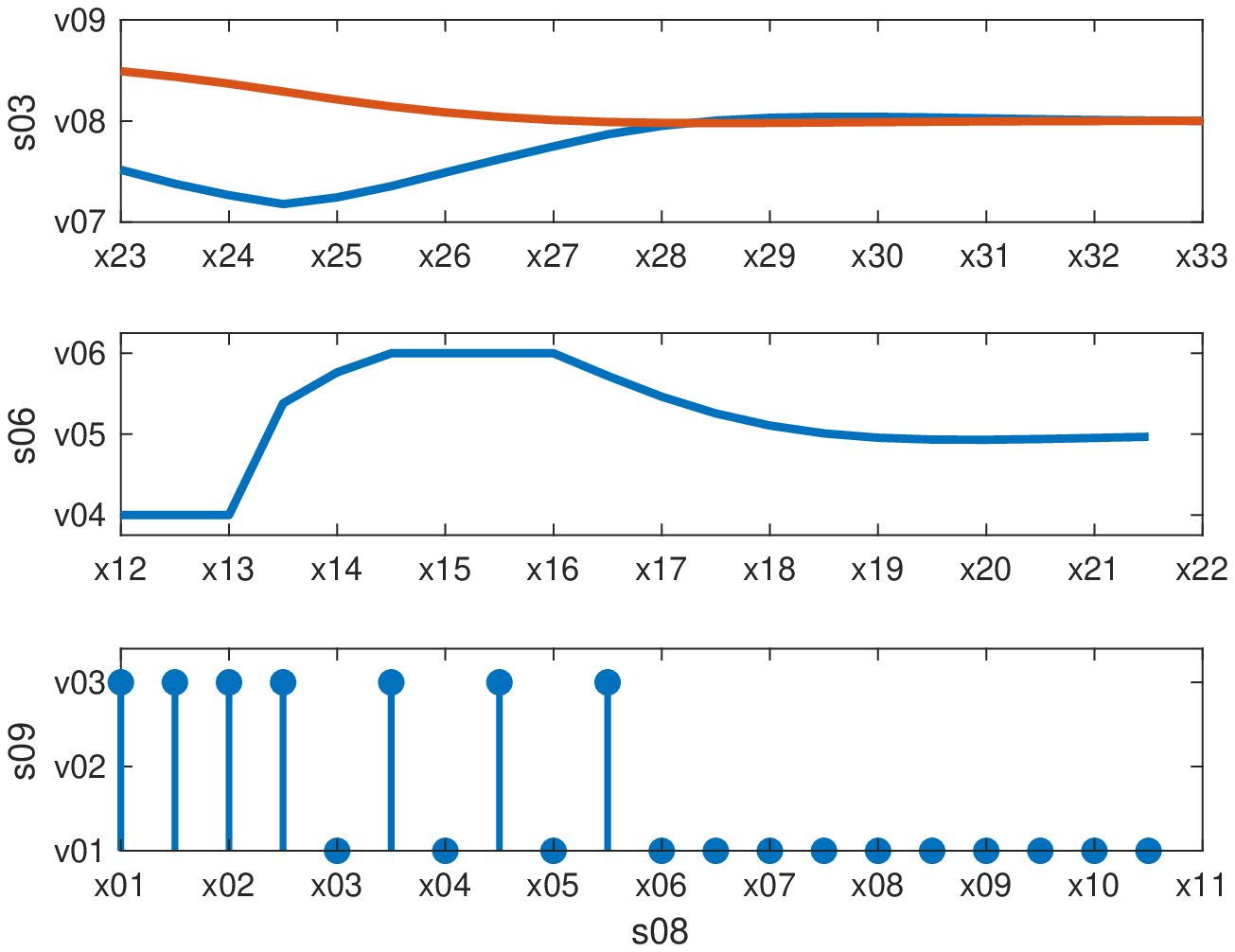}%
\end{psfrags}%
%
% End Seron2003RegA2.tex
}\quad \quad
\subfloat[approach proposed here]{% This file is generated by the MATLAB m-file laprint.m. It can be included
% into LaTeX documents using the packages graphicx, color and psfrag.
% It is accompanied by a postscript file. A sample LaTeX file is:
%    \documentclass{article}\usepackage{graphicx,color,psfrag}
%    \begin{document}\input{Seron2003OnA2}\end{document}
% See http://www.mathworks.de/matlabcentral/fileexchange/loadFile.do?objectId=4638
% for recent versions of laprint.m.
%
% created by:           LaPrint version 3.16 (13.9.2004)
% created on:           22-Oct-2018 08:57:35
% eps bounding box:     15 cm x 11.25 cm
% comment:              
%
\begin{psfrags}%
\psfragscanon%
\scriptsize%
%
% text strings:
\psfrag{s02}[b][b]{\color[rgb]{0.15,0.15,0.15}\setlength{\tabcolsep}{0pt}\begin{tabular}{c}$x(k)$\end{tabular}}%
\psfrag{s05}[b][b]{\color[rgb]{0.15,0.15,0.15}\setlength{\tabcolsep}{0pt}\begin{tabular}{c}$u(k)$\end{tabular}}%
\psfrag{s07}[t][t]{\color[rgb]{0.15,0.15,0.15}\setlength{\tabcolsep}{0pt}\begin{tabular}{c}$k$\end{tabular}}%
\psfrag{s08}[b][b]{\color[rgb]{0.15,0.15,0.15}\setlength{\tabcolsep}{0pt}\begin{tabular}{c}$e(k)$\end{tabular}}%
%
% axes ticklabel color:
\color[rgb]{0.15,0.15,0.15}%
%
% xticklabels:
\psfrag{x01}[t][t]{0}%
\psfrag{x02}[t][t]{2}%
\psfrag{x03}[t][t]{4}%
\psfrag{x04}[t][t]{6}%
\psfrag{x05}[t][t]{8}%
\psfrag{x06}[t][t]{10}%
\psfrag{x07}[t][t]{12}%
\psfrag{x08}[t][t]{14}%
\psfrag{x09}[t][t]{16}%
\psfrag{x10}[t][t]{18}%
\psfrag{x11}[t][t]{20}%
\psfrag{x12}[t][t]{}%
\psfrag{x13}[t][t]{}%
\psfrag{x14}[t][t]{}%
\psfrag{x15}[t][t]{}%
\psfrag{x16}[t][t]{}%
\psfrag{x17}[t][t]{}%
\psfrag{x18}[t][t]{}%
\psfrag{x19}[t][t]{}%
\psfrag{x20}[t][t]{}%
\psfrag{x21}[t][t]{}%
\psfrag{x22}[t][t]{}%
\psfrag{x23}[t][t]{}%
\psfrag{x24}[t][t]{}%
\psfrag{x25}[t][t]{}%
\psfrag{x26}[t][t]{}%
\psfrag{x27}[t][t]{}%
\psfrag{x28}[t][t]{}%
\psfrag{x29}[t][t]{}%
\psfrag{x30}[t][t]{}%
\psfrag{x31}[t][t]{}%
\psfrag{x32}[t][t]{}%
\psfrag{x33}[t][t]{}%
%
% yticklabels:
\psfrag{v01}[r][r]{0}%
\psfrag{v02}[r][r]{}%
\psfrag{v03}[r][r]{1}%
\psfrag{v04}[r][r]{-2}%
\psfrag{v05}[r][r]{0}%
\psfrag{v06}[r][r]{2}%
\psfrag{v07}[r][r]{-2}%
\psfrag{v08}[r][r]{0}%
\psfrag{v09}[r][r]{2}%
%
% Figure:
\includegraphics[width=0.3\textwidth]{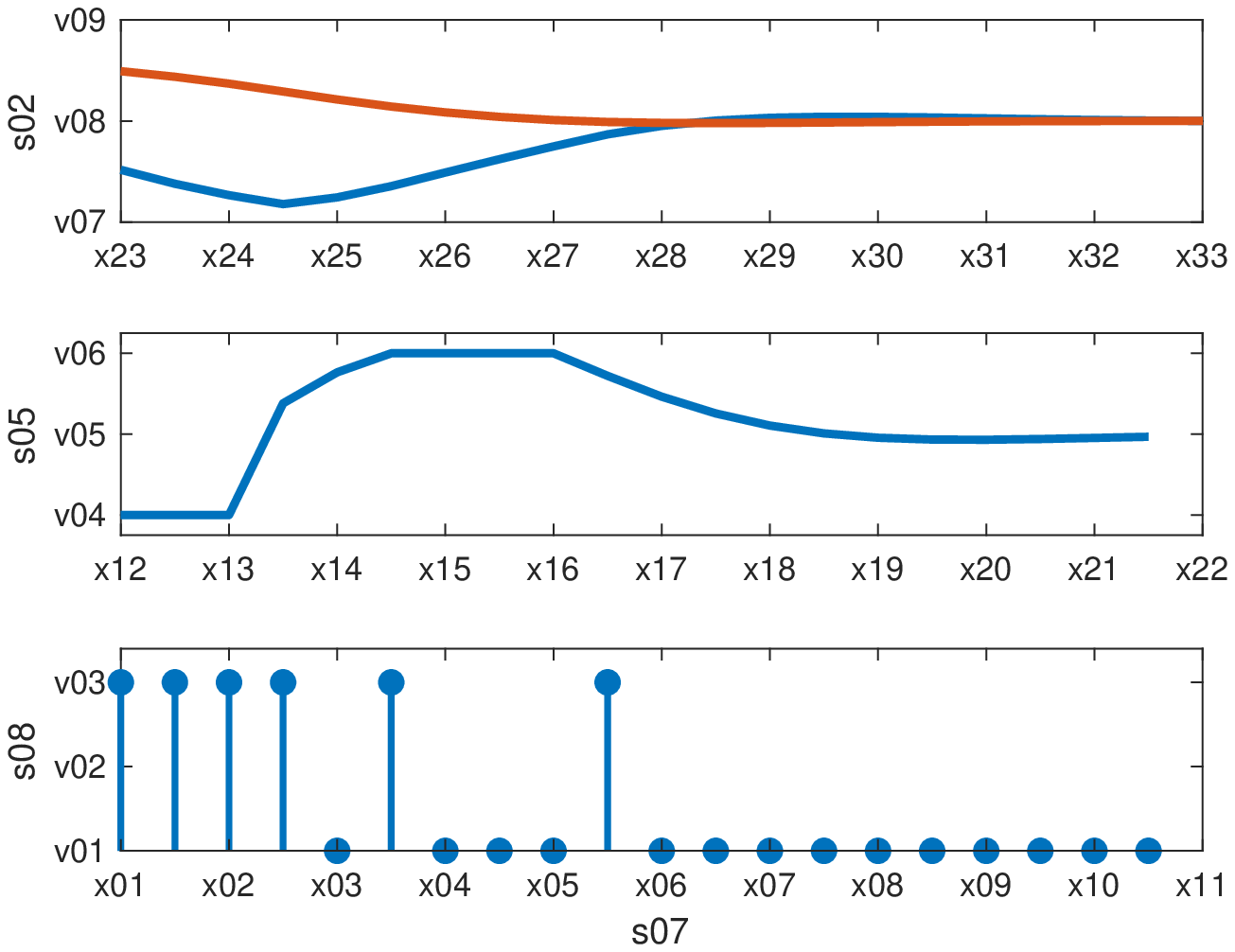}%
\end{psfrags}%
%
% End Seron2003OnA2.tex
}\\
\caption{Results for Example~\ref{example:SISO} for an arbitrary initial state and (a) assuming $\Gamma(\mathcal{A})$ as defined in~\eqref{eq:UnionOfPolytopes} is known for all active sets $\mathcal{A}$ that exist; (b) generated with the existing approach~\cite{Jost2015a}; (c) generated with the approach proposed in Sect.~\ref{sec:Approaches}. The upper plots illustrate the closed-loop trajectory and polytopes in state-space. Blue triangles mark time steps in which a QP must be solved. White circles indicate time steps in which a feedback law can be reused. The bottom plots show the trajectories of the states $x(k)$, inputs $u(k)$ and the function $e(k)$ which indicates whether a QP is solved ($e(k)= 1$) in step $k$ or not ($e(k)= 0$). The results for 1000 random initial states are given in Table \ref{tab:Ergebnisse}. Note that 1536 QPs (\unit[29]{\%}) can be saved with the approach proposed here compared to Jost et al. \cite{Jost2015a}.} 
\label{fig:BeispielSISO}
\end{figure*}

Figure \ref{fig:BeispielSISO} illustrates results for Example~\ref{example:SISO} for an arbitrary initial state. First note that all three approaches result in the same input signal sequence (cp.\ the time series shown in the middle) and closed-loop trajectory (cp.\ the top time series). 
Part (a) of the figure shows the terminal region (cyan polytope), two sets of polytopes $\Gamma(\mathcal{A})$ (red and magenta polytopes), and a single polytope (yellow) 
through which the selected trajectory passes. The polytopes in the magenta and red regions, which each consist of 26 polytopes, have the feedback laws $u= -2$ and $u= 2$ in common, respectively. Consequently, the feedback law does not change in the first three time steps, and it does not change in time steps 5 to 8. This is also evident from the bottom time series shown in Figure~\ref{fig:BeispielSISO} (a), where $e(k)$ indicates a QP is solved in time steps $k= 0$ and $k= 5$ but not in $k= 1, 2$ and $k= 6, 7, 8$. 
%\begin{landscape}
It is instructive to first analyze time steps five to eight in Figures \ref{fig:BeispielSISO} (b) and \ref{fig:BeispielSISO} (c) and compare them to the result in Figure \ref{fig:BeispielSISO} (a):
\begin{itemize}
\item In Figure~\ref{fig:BeispielSISO} (b), i.e., for the approach proposed by Jost et al. \cite{Jost2015a}, the solution of the QP at time step five results in the active set $\mathcal{A}=\{1,7,13\}$ and the feedback law $u=2$. The feedback law can be reused for only one time step because only one polytope is calculated at time step five. Although the feedback law does not change from time step six to eight, a new QP is solved in time step seven to determine the subsequent polytopes. 
\item In Figure~\ref{fig:BeispielSISO} (c), i.e., for the approach proposed here, the feedback law $u=2$ determined in time step five can be reused for three time steps. 
The feedback law is uniquely defined by the active subset $\tilde{\mathcal{A}}=\{1\}$. 
Consequently, the set $\mathcal{E}(\mathcal{A})=\{\{1,7,13\},\{1,7\},\{1,13\},\{1\}\}$ defined by \eqref{eq:SetsForSubregion} yields three non-empty polytopes in time step five. 
As a result, the pointwise solution at time step five is exploited more widely here than in the approach proposed by Jost et al. \cite{Jost2015a} and thus results in a larger number of polytopes for the feedback law.  
\end{itemize}

In contrast to time steps five to eight, no additional savings can be achieved in the first three steps. In fact, the approach proposed here does find several polytopes with the same feedback law as for $k= 0$ (three topmost red polytopes in Figure~\ref{fig:BeispielSISO} (c)), but the trajectory does not pass through them. As a result, the same QPs need to be solved in the approach proposed here as in Jost et al. \cite{Jost2015a}. 

Since results on a single trajectory are only anecdotal, we consider 1000 random initial states. The results are summarized in Table \ref{tab:Ergebnisse}. 
The table states the reusability of feedback laws in percent. Since the feedback law in the terminal set, i.e., the unconstrained linear-quadratic regulator, can be reused arbitrarily long once the terminal set has been entered, we exclude steps in the terminal set in our comparison\footnote{Note that their consideration would improve the statistics in favor of the approach proposed here.}. 
The table shows that in the best case, i.e., assuming $\Gamma(\mathcal{A})$ are known for all existing active sets $\mathcal{A}$, a reusability of \unit[39]{\%} can be achieved. The approach from Jost et al. \cite{Jost2015a} attains a reusability of \unit[7.2]{\%} only. The approach proposed here, in contrast, achieves a reusability of \unit[34]{\%}, which amounts to \unit[87]{\%} of the maximal value reached if all $\Gamma(\mathcal{A})$ were known. 
Note that the proposed approach reduces the number of QPs by 1536 (\unit[29]{\%}) and the computation time by about \unit{23}{\%} compared to Jost et al. \cite{Jost2015a}.\footnote{Computation times are matlab execution times.}  

Table \ref{tab:Ergebnisse} also gives the results for Example~\ref{example:MIMO}. 
If all active sets $\mathcal{A}$ that exist for the problem and $\Gamma(\mathcal{A})$ were determined beforehand, a reusability of \unit[12]{\%} could be achieved. 
We stress again the proposed approach does not require these offline calculations, but they are carried out only for the sake of a comparison only. In the approach from Jost et al. \cite{Jost2015a}, no feedback law at all is reused for the 1000 random initial conditions. In the approach proposed here, the optimal feedback law is reused in \unit[11]{\%} of all cases, which amounts to \unit[89]{\%} of the maximal achievable reusability. The computation time can be reduced by about \unit{6.3}{\%} with proposed approach compared to Jost et al. \cite{Jost2015a}. Moreover, active sets arise that satisfy the conditions in Proposition \ref{prop:uUniquelyDetermined} for non-saturated feedback laws.

The results for Example~\ref{example:InvertedPendulum} are summarized in Table \ref{tab:Ergebnisse}. 
If all active sets $\mathcal{A}$ and all sets $\Gamma(\mathcal{A})$ were calculated explicitly, 
a reusability of $\unit[41]{\%}$ could be achieved. The approach proposed here, in contrast, reuses feedback laws and avoids solving QPs in \unit{41}{\%} of the cases without computing and storing the explicit solution or parts thereof. Compared to Jost et al. \cite{Jost2015a}, which fails to reuse feedback laws for this example, the number of QPs can be reduced by 1434 and the computation time by about \unit{23}{\%}.
\setlength{\tabcolsep}{2pt}
\begin{table}[]
\scriptsize
\begin{center}
\caption{Maximum reusability of feedback laws (2nd column) and reusabilities achieved with the approach from Jost et al. \cite{Jost2015a} (3rd column) and the approach proposed in Section \ref{sec:Approaches} (4th column) for 1000 random initial states each for examples 1 and 2.  }
\begin{tabular}{|c|c|c|c|} \hline
     & \ \ $\Gamma(\mathcal{A})$ known \ \ & \multicolumn{2}{|c|}{online approaches}  \\ \cline{3-4} 
        example  & \ \ for all ex.\ $\mathcal{A}$ \ \ & \ \ \cite{Jost2015a} \ \ &\ \ proposed here \\ \hline
        1 & \unit[39]{\%} &  \unit[7.2]{\%} & \unit[34]{\%} \\ \hline
         2 & \unit[12]{\%} &  \unit[0]{\%} &  \unit[11]{\%} \\ \hline    
         3 & \unit[41]{\%} &  \unit[0]{\%} & \unit[41]{\%} \\ \hline        
\end{tabular}
\label{tab:Ergebnisse}
\end{center}
\end{table}
We point out that polytopes may be calculated in the proposed approach that are not passed by the closed-loop trajectory (see five out of eight red polytopes in Figure \ref{fig:BeispielSISO} (c), for example). Since solving a QP~\eqref{eq:QP} is much more expensive than computing a polytope with \eqref{eq:Td}, the computation of a small number of additional polytopes is acceptable in that overall computational savings result. 

	We implement the approach proposed in Section \ref{sec:Approaches} in a networked MPC variant~\cite{BernerP2016}
	to analyze its usefulness for embedded hardware. 
	QPs are solved on a computationally powerful central node on demand in this setting. 
	Active sets are transmitted to a lean local node, 
	where input signals are computed by evaluating the optimal affine feedback law. 
	Whenever the current affine law is not optimal anymore, the central node is requested to solve a new QP. 

	In contrast to Berner and M\"onnigmann \cite{BernerP2016}, the central node does not compute and transmit a single active set here,
	but the set $\mathcal{E}(\mathcal{A})$ or a subset thereof. 
	Apart from the reduction in the number of QPs that need to be solved, this
	reduces the bandwidth requirements of the networked MPC more efficiently than in Berner and M\"onnigmann \cite{BernerP2016}. 
	We can control the amount of transmitted data 
	by limiting the number of active sets in $\mathcal{E}(\mathcal{A})$. 
	In our implementation an active set is represented as a tuple of $q$ bits $\alpha=(\alpha_q, \ldots, \alpha_1)$, 
	where $\alpha_i=1$ if $i \in \mathcal{A}$ and $\alpha_i=0$ otherwise. 
	The active sets in $\mathcal{E}(\mathcal{A})$ are sorted by descending binary number of their tuple 
	and only the first $l$ sets are transmitted to the local node.

As the central node, we use a standard desktop computer with Intel Core i5-8400 CPU with 2.8 GHz and 8GB RAM. The central node is connected to an IEEE 802.11 b/g/n wireless LAN access point. The local node is an Espressif ESP8266 SoC with an integrated IEEE 802.11 b/g/n WiFi controller. The SoC features a 80MHz Tensilica L106 32-bit RISC micro controller and 96 KiB data RAM. In our implementation the local node generates closed-loop control signals for the inverted pendulum example for 1000 random initial states. 
If we set the maximum number of transmitted sets per request to $l=50$, the number of requests to the central node can be reduced by about \unit[37.6]{\%} compared to Jost et al. \cite{Jost2015a}. If a maximum number of $l=10$ sets is chosen, the reduction in the number of requests is about \unit[33.1]{\%}. It is remarkable that with a maximum number of $l=5$ sets there is still a reduction in the number of requests of about \unit[27.5]{\%}.

\section{Conclusions}\label{sec:conclusion}
	We introduced a simple criterion for finding, from the solution of the MPC problem for the current state, 
	state-space regions with the same optimal affine $x\rightarrow u(x)$ as the current state. 
	Since the function $x\rightarrow u(x)$ can be determined from the solution for the current state 
	with very small additional computational effort, 
	the MPC solution can be determined on a full-dimensional state-space region from a single QP,
	where multiple QPs would have to be solved otherwise. 

	In contrast to existing methods we neither reuse $x\rightarrow u(x)$ on only a single polytope, 
	nor reuse it as an approximately optimal feedback on neighbouring polytopes. 
	Essentially, we exploited that the solution to the MPC problem is often defined 
	by the activity of the same subset of the constraints on a union of polytopes. 
	We proposed a heuristics that finds a subset of this union of polytopes at runtime and, by reusing the optimal feedback law as long as the system stays in this union, reduces the computational effort of MPC by avoiding obsolete QPs. We stress the heuristics exploits the piecewise-affine character without requiring the explicit solution. 
	The reported computational experiments showed that computational effort can be reduced considerably
	compared to an existing approach.

Reducing the number of QPs particularly has a benefit in a networked setting as proposed in Berner and M\"onnigmann \cite{BernerP2016}, because it reduces the number of requests to a central node. 

\subsection*{Acknowledgement}
Support by the Deutsche Forschungsgemeinschaft (DFG) under grant MO 1086/15-1 is gratefully acknowledged.

\bibliographystyle{plain}
\bibliography{literature}%

\end{document}